\documentclass [reqno]{amsart}
\usepackage{amsmath}
\usepackage{amsfonts,amssymb,amscd,amsmath, graphicx, bbm}
\usepackage{url}

\newtheorem{remark}{Remark}

\begin{document}

\title [Series involving central binomial coefficients \& harmonic numbers]
{\bfseries\large{Sums of series involving central binomial coefficients \& harmonic numbers}}

\date{\today}

\author {Amrik Singh Nimbran}
\address{B3-304, Palm Grove Heights, Ardee City, Gurgaon, Haryana, INDIA }
\email{amrikn622@gmail.com}


\begin{abstract}
This paper contains a number of series whose coefficients are products of central binomial coefficients \& harmonic numbers. An elegant sum involving $\zeta(2)$ and two other nice sums appear in the last section.
\end{abstract}

\maketitle

\section{Introduction and Preliminary resullts}

\subsection{Beginnings}

Euler investigated the partial sums of the harmonic series, finding connection between them and $\log (n).$ These sums: $\displaystyle H_n=\sum_{k=1}^{n} \frac{1}{k}=\int_{0}^1 \frac{1-x^n}{1-x} \, dx$ are generally known as \emph{harmonic numbers.} He then discovered the formula for $\displaystyle \zeta(2m)=\sum_{n=1}^\infty \frac{1}{n^{2m}}.$ He also introduced the dilogarithm function:
\[
{\rm Li}_2 (x)=\sum_{m=1}^\infty \frac{x^m}{m^2}=-\int_{0}^x \frac{\log(1-t)}{t} {dt}=-\int_{0}^x \frac{dt}{t} \int_{0}^t \frac{du}{1-u} \qquad (|x|\le1)
\]
and computed these values: $\displaystyle
{\rm Li}_2 (1) = \frac{\pi^2}{6} , \; {\rm Li}_2 (-1) = -\frac{\pi^2}{12}, \; {\rm Li}_2 \left(\frac{1}{2}\right) =\frac{\pi^2}{12} -\frac{\log^2 2}{2}.$

Further, he investigated the double sums $\displaystyle \sum_{n=1}^\infty \frac{H_n}{n^m} \quad (\text {for some fixed}\;  m\ge 2)$ which have been a popular topic of research in recent years.

\subsection{Generating functions for central binomial coefficients and Catalan numbers}

Recall that the \emph{generating function} for the sequence $a_0, a_1, a_2, \dots$ is defined to be the function represented by power series: $\displaystyle G(x):=\sum_{n=0}^\infty a_n \, x^n.$ It is always permissible to integrate (and differentiate) a power series term by term over any closed interval lying entirely within its interval of convergence.

The binomial coefficient $\displaystyle \binom{2n}{n}$, the largest coefficient of the polynomial $(1 +x)^{2n},$ forms the central column of Pascal's triangle and so is always an integer. The sequence of these numbers is generated by
\begin{equation}
\frac{1}{\sqrt{1-4x}}=\sum_{n=0}^\infty \binom{2n}{n} x^n. \label{eq:1}
\end{equation}
The expansion is a consequence of the binomial theorem as for  $n \in \mathbb{N}:$
\[
\frac{1}{(1 -x)^{1/n}} =1 +\frac{1}{n} x +\frac{1\cdot(1+n)}{n\cdot2n} x^2 +\frac{1\cdot(1+n)\cdot(1+2n)}{n\cdot2n\cdot3n} x^3 +\dots.
\]

The series on the R.H.S. of (\ref{eq:1}) converges if $|x|<\frac{1}{4}.$ Lehmer\cite{lehmer} obtained some `interesting series' by repeated integrations of (\ref{eq:1})

If we integrate (\ref{eq:1}) from 0 to $x$ nd then divide the result by $x$ we get generating function for $\frac{1}{n+1}\binom{2n}{n},$ known as the \emph{Catalan numbers} and denoted by $C_n.$
\begin{equation}
\frac{1-\sqrt{1-4x}}{2x}=\sum_{n=0}^\infty C_n x^n. \label{eq:1a}
\end{equation}

Transposing the first term of the right side of (\ref{eq:1}) to the left, dividing both sides by $x$ and then integrating, one gets \cite[(6)]{lehmer}:
\begin{equation}
\sum_{n=1}^\infty \frac{1}{n}\binom{2n}{n} x^n =2 \log \frac{1-\sqrt{1-4x}}{2x}. \label{eq:1b}
\end{equation}

\subsection{Generating function for harmonic numbers and few series}
Using Euler's integral representation of harmonic numbers, we have:
\begin{align*}
\sum_{n=1}^\infty H_n \, x^n \, {dx}& =\sum_{n=1}^\infty \left(\int_{0}^1 \frac{1-u^n}{1-u} \, {du}\right) x^n =\int_{0}^1 \frac{1}{1-u} \left(\sum_{n=1}^\infty (x^n -(ux)^n) \right) {du}\\
& =\int_{0}^1 \frac{1}{1-u} \left(\frac{1}{1-x} -\frac{1}{1-ux} \right) {du} =x\int_{0}^1 \frac{du}{1-ux}.
\end{align*}
Since
\[
\frac{x}{1-x}\int_{0}^1 \frac{du}{1-ux} =-\frac{\ln (1-x)}{1-x},
\]
we get the generating function for harmonic numbers given in \cite[p.54, 1.514.6]{gradshteyn}:
\begin{equation}
-\frac{\ln (1-x)}{1-x} =\sum_{n=1}^\infty H_n x^n \qquad  (x^2<1).\label{eq:1c}
\end{equation}

Considering the partial sums of the alternating harmonic series (having sum $\log 2$), with notation
$\displaystyle H_n^{\prime}=\sum_{k=1}^n \frac{(-1)^{k+1}}{k}=\log 2 +\frac{(-1)^n}{2}\left[\psi\left(\frac{n+1}{2}\right) -\psi\left(\frac{n+2}{2}\right)\right]$ where $\psi(x)$ is the digamma function, we have this generating function:
\begin{equation}
\frac{\ln (1+x)}{1-x} =\sum_{n=1}^\infty H_n^{\prime}\,  x^n \qquad  (x\ne 1).\label{eq:1d}
\end{equation}

Now
\[
\int \frac{\log(1 + x)}{1 - x} \, {dx} =-{\rm Li}_2\left(\frac{1+x}{2}\right) -\log\left(\frac{1 - x}{2}\right) \log(1+x) +C.
\]
This integral evaluated between -1 and 0 gives the sum:
\begin{equation}
\sum_{n=1}^\infty (-1)^{n+1} \frac{H_n^{\prime}}{n+1} =\frac{\pi^2}{12} -\frac{1}{2} \log^{2} 2. \label{eq:1e}
\end{equation}

Further,
\begin{align*}
\int \frac{\log(1 + x)}{x(1 -x)} \, {dx} &=\int \frac{\log(1 + x)}{1 -x} \, {dx}+\int \frac{\log(1 + x)}{x} \, {dx}\\
&=-{\rm Li}_2(x) -{\rm Li}_2\left(\frac{1+x}{2}\right) -\log\left(\frac{1 -x}{2}\right) \log(1+x) +C
\end{align*}
which yields the sum valid for $-1\le x<1:$
\begin{equation}
\sum_{n=1}^\infty \frac{H_n^{\prime}}{n} \, x^n ={\rm Li}_2\left(\frac{1}{2}\right)-{\rm Li}_2(-x) -{\rm Li}_2\left(\frac{1+x}{2}\right)
-\log\left(\frac{1-x}{2}\right) \log (1+x). \label{eq:1f}
\end{equation}
Its alternative form occurs in \cite[p.302, A.2.8 (1)]{lewin} \cite[(5)]{boyadzhiev13}
\[
\sum_{n=1}^\infty \frac{H_n^{\prime}}{n}\, x^n={\rm Li}_2\left(\frac{1-x}{2}\right)-{\rm Li}_2\left(\frac{1}{2}\right)-{\rm Li}_2(-x) -\log(1-x) \log 2
\]
and Boyadzhiev \cite[(10)]{boyadzhiev13} computes the sum differently:
\begin{equation}
\sum_{n=1}^\infty (-1)^{n+1} \frac{H_n^{\prime}}{n} =\frac{\pi^2}{12} +\frac{1}{2} \log^{2} 2. \label{eq:1g}
\end{equation}

So the sum and difference of (\ref{eq:1e}) and (\ref{eq:1g}) result in
\[
\sum_{n=1}^\infty (-1)^{n+1} \frac{H_n^{\prime}}{n(n+1)} = \log^{2} 2; \quad \sum_{n=1}^\infty (-1)^{n+1} \frac{(2n+1) \,H_n^{\prime}}{n(n+1)} =\frac{\pi^2}{6}.
\]

The two integrals preceding  equations (\ref{eq:1e}) and (\ref{eq:1f}) evaluated between 0 and $1/2$ yield sums:
\begin{equation}
\sum_{n=1}^\infty  \frac{H_n^{\prime}}{n \, 2^{n}}=\frac{1}{2}{\rm Li}_2\left(\frac{1}{4}\right) +\log^{2} 2
=-\frac{1}{6}{\rm Li}_2\left(\frac{1}{9}\right) +\frac{\pi^2}{36} -\frac{\log^{2} 3}{3} +\log 2 \log 3, \label{eq:1h}
\end{equation}
and
\begin{equation}
\sum_{n=1}^\infty  \frac{H_n^{\prime}}{(n+1) \, 2^{n+1}}=-\frac{1}{3}{\rm Li}_2\left(\frac{1}{9}\right) -\frac{\pi^2}{36}  -\frac{1}{2} \log^{2} 2 -\frac{2}{3} \log^2 3 +2 \log 2 \log 3. \label{eq:1i}
\end{equation}

Multiplying (\ref{eq:1i}) by 2 and by subtracting the result from (\ref{eq:1h}), we get
\begin{equation}
\sum_{n=1}^\infty  \frac{(3n+4) H_n^{\prime}}{n(n+1) \, 2^{n}}=\frac{2\pi^2}{3} +(\log 2)^2=8\sum_{n=1}^\infty (-1)^{n+1} \frac{H_n^{\prime}}{n}. \label{eq:1j}
\end{equation}

\subsection{Relation between binomial coefficients and harmonic numbers}

Harmonic numbers can be expressed in term of binomial coefficients:
\begin{align*}
H_n &=\sum_{k=0}^{n-1} \frac{1}{k+1}=\sum_{k=0}^{n-1} \left(\int_{0}^1 x^{k} \, {dx}\right) =\int_{0}^1 \left(\sum_{k=0}^{n-1} x^k \right) \, {dx}\\
&=\int_{0}^1 \frac{1 -x^n}{1-x} \, {dx} \qquad (\text{partial sum of the G.P.})\\
&=\int_{0}^1 \frac{1 -(1-y)^n}{y} \, {dy} \qquad (\text{on setting} \; 1-x=y)\\
&=\int_{0}^1 \frac{\sum_{k=1}^n \binom{n}{k} \, (-1)^{k-1} \, y^k}{y} \, {dy} \qquad (\text{removing 1 from numerator})\\
&=\sum_{k=1}^n \binom{n}{k} (-1)^{k-1} \int_{0}^1 y^{k-1} \, {dy}=\sum_{k=1}^n \binom{n}{k} (-1)^{k-1} \frac{1}{k}.
\end{align*}

No $H_n$ is an integer for $n>1.$ We find  in \cite[p.192, (5.48)]{graham} this inversion formula:
\[
g(n)=\sum_{k} \binom{n}{k} (-1)^k f(k) \Longleftrightarrow f(n)=\sum_{k} \binom{n}{k} (-1)^k g(k).
\]

Since $\displaystyle H_n=\sum_{k=1}^n (-1)^k \binom{n}{k} \frac{1}{k},$ we have: $\displaystyle \frac{1}{n}=\sum_{k=1}^n (-1)^k \frac{1}{k} \binom{n}{k} H_k.$

This relation is used by Boyadzhiev \cite{boyadzhiev11} \cite{boyadzhiev12} for derivation of certain series whose coefficients are products of the central binomial coefficients and harmonic numbers, and whose sums are expressible in term of logarithms.

We intend to obtain here series whose sums involve $\pi, \, \zeta(2)$ and Catalan's constant and thereby supplement Boyadzhiev's work.

\section{Generating function for $\binom{2n}{n} \, H_n$ and known series}

\subsection{Using Euler's transformation of series}
We find in \cite[p.469]{knopp} this version of Euler's transformation of series:
\[
\sum_{k=0}^\infty a_k x^{k+1}=\sum_{k=0}^\infty a_k \left (\frac{y}{1-y} \right)^{k+1}=\sum_{n=0}^\infty \left\{\sum_{m=0}^n \binom{n}{m} a_m \right\} y^{n+1}.
\]

Boyadzhiev comes up with this formula for $\alpha \in \mathbb{C}$ in \cite[(2.4)]{boyadzhiev11} and \cite[(10)]{boyadzhiev12}:

\[
\sum_{n=0}^\infty \binom{\alpha}{n} (-1)^{n} a_n \, z^{n} =(z+1)^{\alpha} \sum_{n=0}^\infty \binom{\alpha}{n} (-1)^{n} \left\{\sum_{k=0}^n \binom{n}{k} a_k \right\} \, \left(\frac{z}{z+1}\right)^{n}
\]

Setting $z=4x, \, a_k=(-1)^{k-1} H_k, \, \alpha=-1/2$ and using the relation we derived above, Boyadzhiev obtained for $|x|<\frac{1}{4}$ these generating functions for the product of the harmonic numbers and the central binomial coefficients:
\begin{gather}
\sum_{n=0}^\infty H_n \, \binom{2n}{n} \, (-1)^{n+1} x^{n}=\frac{2}{\sqrt{1+4x}} \log \left(\frac{2\sqrt{1+4x}}
{1+\sqrt{1+4x}} \right). \label{eq:2}\\
\sum_{n=0}^\infty H_n \, \binom{2n}{n} \, x^{n}=\frac{2}{\sqrt{1-4x}} \log \left(\frac{1+\sqrt{1-4x}}
{2\sqrt{1-4x}} \right). \label{eq:3}
\end{gather}

\subsection{Alternative derivation}

We can also derive (\ref{eq:3}) in a different way as follows:
\begin{align*}
\sum_{n=1}^\infty \binom{2n}{n} \, H_n \, x^{n} &=\sum_{n=1}^\infty \binom{2n}{n} \, x^{n} \,\int_{0}^1 \frac{1-t^n}{1-t} \, {dt}\\
& =\int_{0}^1 \frac{1}{1-t} \, \left(\sum _{n=1}^{\infty} \binom{2n}{n} (x^n -(xt)^n)  \right) {dt}\\
& =\int_{0}^1 \frac{1}{1-t} \, \left(\sum _{n=1}^{\infty} \binom{2n}{n} x^n  -\sum _{n=1}^{\infty} \binom{2n}{n} (xt)^n \right) {dt} \\
& =\int_{0}^1 \frac{1}{1-t} \, \left(\frac{1}{\sqrt{1-4x}} -\frac{1}{\sqrt{1-4xt}}\right) {dt}  \qquad [\text{by using (\ref{eq:1})}].
\end{align*}
We assumed that swapping of summation and integration is permissible here.

Now $\displaystyle \frac{1}{\sqrt{1-4x}} \int \frac{dt}{1-t}=-\frac{\log(1-t)}{\sqrt{1-4x}}.$
To evaluate the second integral we make the substitution: $u=\sqrt{1-4xt}$ so that $\displaystyle du=-\frac{2x}{\sqrt{1-4xt}} {dt}$ and
$\displaystyle 1-t=\frac{4x-1+u^2}{4x}.$ Then,
\[
\int \frac{dt}{(1-t)\sqrt{1-4xt}}=2x\int \frac{du}{(1-4x) -u^2}=\frac{1}{\sqrt{1-4x}} \log \frac{\sqrt{1-4x} +\sqrt{1-4xt}}{\sqrt{1-4x} -\sqrt{1-4xt}}.
\]
Thus the integral becomes:
\begin{align*}
&-\frac{1}{\sqrt{1-4x}} \log \frac{(\sqrt{1-4x} +\sqrt{1-4xt})(1-t)}{(\sqrt{1-4x} -\sqrt{1-4xt})}\\
&=-\frac{2}{\sqrt{1-4x}} \log \frac{(\sqrt{1-4x} +\sqrt{1-4xt})(1-t)}{-4x(1-t)}
\end{align*}
that is,
\[
-\frac{2}{\sqrt{1-4x}} \log \frac{\sqrt{1-4x} +\sqrt{1-4xt}}{-4x}.
\]
which at $t=1$ has the value $\displaystyle -\frac{2}{\sqrt{1-4x}} \log \frac{\sqrt{1-4x}}{-2x}$ and at $t=0$ it becomes:
$\displaystyle \frac{2}{\sqrt{1-4x}} \log \frac{1+\sqrt{1-4x}}{-4x}.$ Thus the definite integral becomes:
\[
\frac{2}{\sqrt{1-4x}} \log \frac{1+\sqrt{1-4x}}{-4x}  \frac{-2x}{\sqrt{1-4x}}
=\frac{2}{\sqrt{1-4x}} \log \frac{1+\sqrt{1-4x}}{2\sqrt{1-4x}}
\]
which is the formula (\ref{eq:3}). $\Box$

Integrating the power series (\ref{eq:3}), using the substitution $1-4x=y^{2}$ for the RHS, one obtains for every $|x|\le \frac{1}{4}$,
\begin{align}
\sum _{n=0}^{\infty} H_n \, \binom{2n}{n}\frac{x^{n+1}}{n+1} & =\sqrt{1-4x} \log(2\sqrt{1-4x}) \nonumber\\
&\quad -(1+\sqrt{1-4x}) \log(1+\sqrt{1-4x})+ \log 2. \label{eq:4}
\end{align}

Putting $x=1/4$ in (\ref{eq:4}) yields:
\begin{equation}
\sum_{n=1}^\infty \frac{H_n \, \binom{2n}{n}}{(n+1) \, 2^{2n}} =4 \log 2. \label{eq:5}
\end{equation}

By shifting the index we get:
\begin{align*}
\sum_{n=1}^\infty \frac{H_n \, \binom{2n}{n}}{(2n-1) \, 2^{2n}} &=\sum_{n=0}^\infty \frac{H_{n+1} \, \binom{2n+2}{n+1}}{(2n+1) \, 2^{2n+2}}
=\frac{1}{2}\left[\sum_{n=0}^\infty \frac{(H_{n} +\frac{1}{n+1}) \,\binom{2n}{n}}{(n+1)\, 2^{2n}}\right]\\
&=\frac{1}{2}\left[\sum_{n=0}^\infty \frac{H_{n} \,\binom{2n}{n}}{(n+1)\, 2^{2n}} +\sum_{n=0}^\infty \frac{\,\binom{2n}{n}}{(n+1)^2\, 2^{2n}}\right]
\end{align*}
and since (see \cite[pp.251-252, \S 1081]{edwards})
\begin{equation}
\sum_{n=0}^\infty \frac{\,\binom{2n}{n}}{(n+1)^2\, 2^{2n}}=4 -4\log 2 \label{eq:6}
\end{equation}
we obtain by using (\ref{eq:5}) and (\ref{eq:6}):
\begin{equation}
\sum_{n=1}^\infty \frac{H_n \, \binom{2n}{n}}{(2n-1) \, 2^{2n}} =2. \label{eq:7}
\end{equation}

Edwards deduced (\ref{eq:6}) via the integral $\displaystyle \int_{0}^1 \frac{\log (1/x)}{\sqrt{1-x}} \, {dx}$ and putting $x=\sin^2\theta.$ He also gave (pp.252--253):$\displaystyle \sum_{n=0}^\infty \frac{\,\binom{2n}{n}}{(2n+1)^2\, 2^{2n}}=\frac{\pi \ln 2}{2}.$ It may be interesting to mention here that $\displaystyle \sum_{n=0}^\infty \frac{\,\binom{2n}{n}}{(2n+1)^2\, 2^{3n}}=\sqrt{2} \left[\frac{\pi \ln 2}{8} +\frac{G}{2}\right],$ where $G$ is Catalan's constant defined by $\displaystyle \sum_{n=0}^\infty \frac{(-1)^n}{(2n+1)^2}\approx 0.9159655941\dots$

We can similarly derive with $x=-\frac{1}{16}$ in (\ref{eq:4}) and by shifting the index:
\begin{equation}
\sum_{n=1}^\infty (-1)^{n+1} \frac{H_n \, \binom{2n}{n}}{(n+1) \, 2^{4n}}=16\log (4\sqrt{5} -8) + 8\sqrt{5} \log (10-4\sqrt{5}), \label{eq:8}
\end{equation}
and
\begin{equation}
\sum_{n=1}^\infty (-1)^{n+1} \frac{H_n \, \binom{2n}{n}}{(2n-1) \, 2^{4n}}=(\sqrt{5}\, -2) +\sqrt{5} \log \left(\frac{1}{\sqrt{5}} +\frac{1}{2}\right). \label{eq:9}
\end{equation}

By a variation on the method, we will now derive some interesting series whose sums involve $\pi, \, \zeta(2)$ and Catalan's constant.

\section{New sums involving $\pi, \zeta (2), G$}

This is the series that we aim to obtain:

\begin{equation}
\sum_{n=1}^\infty \frac{H_n \, \binom{2n}{n}}{n \, 2^{2n}} =\frac{\pi^2}{3}. \label{eq:10}
\end{equation}

Dividing both sides of (\ref{eq:3}) by $x$ and integrating between limits $x=0$ to $x=1/4$, after using the substitution $1-4x=y^{2}, \, x=\frac{1-y^2}{4}, \, {dx}=\frac{-y}{2} \, {dy}$ we get on the R.H.S.:
\[
 I=\int_{1}^{0} \frac{2\cdot 4}{(1-y^2) y} \log \left(\frac{1+y} {2y} \right) \cdot \frac{-y}{2} \, {dy}
=\int_{0}^{1} \frac{4}{(1-y^2)} \log \left(\frac{1+y} {2y} \right) \, {dy}.
\]

The software \emph{WolframAlpha} at \url{https://www.wolframalpha.com} gives:
\begin{align*}
\int 4\log \left(\frac{1+y}{2y}\right) \, \frac{dy}{1-y^2} & =-2{\rm Li}_2 (1-y) -2{\rm Li}_2 (-y) -2{\rm Li}_2 \left(\frac{y+1}{2}\right) \\
& \quad -\log^2(y+1) +2 \log \left(1+\frac{1}{y}\right) \log (1-y) \\
& \quad +2 \log (2) \log(1-y) -2 \log (1-y) \log (y) +C,
\end{align*}
and
\[
 \int_{0}^1 4\log \left(\frac{1+y}{2y}\right) \, \frac{dy}{1-y^2}  =\frac{\pi^2}{3}.
\]
Consequently,
\[
\sum_{n=1}^\infty \frac{H_n \, \binom{2n}{n}}{n \, 2^{2n}} =\frac{\pi^2}{3}.
\]

Let us try to evaluate the integral $I$ by partial fractions decomposition:
\begin{align*}
I & =\int_{0}^{1} \frac{2}{(1-y)} \log \left(\frac{1+y} {2y} \right) \, {dy} \; +\int_{0}^{1} \frac{2}{(1+y)} \log \left(\frac{1+y} {2y} \right) \, {dy}\\
& =2\int_{0}^{1} \frac{\log (1+y)}{1-y} \, {dy}
-2\int_{0}^{1} \frac{\log 2}{1-y} \, {dy} -2\int_{0}^{1} \frac{\log y}{1-y} \, {dy}\\
& +2\int_{0}^{1} \frac{\log (1+y)}{1+y} \, {dy} -2\int_{0}^{1} \frac{\log 2}{1+y} \, {dy} -2\int_{0}^{1} \frac{\log y}{1+y} \, {dy}.
\end{align*}

These are the relevant indefinite integrals:
\begin{align*}
I_1&=\int \frac{\log (1+y)}{1-y} \, {dy}=-{\rm Li}_2 \left(\frac{1+y}{2}\right) -\log \left(\frac{1-y}{2}\right) \, \log (1+y) +C\\
&=-{\rm Li}_2 \left(\frac{1+y}{2}\right) -\log (1-y) \log (1+y) +\log 2 \log (1+y) +C,
\end{align*}

\[
I_2=\int \frac{\log (1+y)}{1+y} \, {dy}=\frac {\log^2 (1+y)}{2} +C,
\]

\[
I_3=\int \frac{\log 2}{1-y} \, {dy}=-\log 2 \, \log (1-y) +C,
\]

\[
I_4=\int \frac{\log 2}{1+y} \, {dy}=\log 2 \, \log (1+y) +C,
\]

\[
I_5=\int \frac{\log y}{1-y} \, {dy}={\rm Li}_2 (1-y) +C,
\]
and
\[
I_6=\int \frac{\log y}{1+y} \, {dy}={\rm Li}_2 (-y) +\log (y) \log (1+y) +C.
\]

\begin{remark}
There arises a problem when one puts $y=1$ in the two integrals (with opposite signs) namely $I_1$ and $I_3$. Though we get two indeterminate terms $\log (0) \times \log (2)$ with opposite signs, we cannot simply cancel them to get 0. Again, when we put $y=0$ in $I_6,$ we get an indeterminate term $\log(0 )\log (1) =\infty \times 0$ which cannot be straightway taken to be 0; for this, we will take limit as $y \to 0.$
\end{remark}

For $I_1 -I_3,$ we notice as in \cite[(6)]{boyadzhiev13}
\[
\frac{d}{dx} {\rm Li_2 \left(\frac{1-x}{2}\right)} = \frac{1}{1-x} \log \frac{1+x}{2} =\frac{\log (1+x)}{1-x} -\frac{\log (2)}{1-x}.
\]
Therefore,
\[
\int_{0}^1 \left(\frac{\log (1+x)}{1-x} -\frac{\log 2}{1-x}\right)={\rm Li_2 \left(\frac{1-x}{2}\right)}\Big\vert_{0}^1=\frac{(\log 2)^2}{2} -\frac{\pi^2}{12}.
\]

For $I_6$ we have: $\displaystyle \lim_{t\to 0} \log t \log (1+t)$ and by applying l'H\^{o}pital's rule, we write:
\begin{align*}
\lim_{t\to 0} \frac{\log t}{1/\log (1-t)}& =\lim_{t\to 0} \frac{(1-t) \log^2(1-t)}{t}\\
&=\lim_{t\to 0} \left[-\log^2(1-t) -2 \log (1-t)\right]=0.
\end{align*}
Thus the indeterminate expression has value 0.

Consequently, we get (on using the values of the logarithm given by Euler) the same value of the integral $I=\frac{\pi^2}{3}$ as returned by the \emph{Wolfram} software. $\Box$

Combining (\ref{eq:5}) and (\ref{eq:10}), we get:
\begin{equation}
\sum_{n=1}^\infty \frac{H_n \, \binom{2n}{n}}{n (n+1)\, 2^{2n}} =\frac{\pi^2}{3} -4 \log 2. \label{eq:11}
\end{equation}

By a shift of the index, we have:
\[
\sum_{n=1}^\infty \frac{H_n \, \binom{2n}{n}}{n \, 2^{2n}}=\sum_{n=0}^\infty \frac{H_{n+1} \, \binom{2n+2}{n+1}}{(n+1) \, 2^{2n+2}}
=\frac{1}{2}\left[\sum_{n=0}^\infty \frac{(H_{n} +\frac{1}{n+1}) \, (2n+1)\binom{2n}{n}}{(n+1)^2 \, 2^{2n}}\right]
\]
and the expression within brackets on the extreme right can be expanded as
\[
2\sum_{n=0}^\infty \frac{H_{n}\binom{2n}{n}}{(n+1)\, 2^{2n}} -\sum_{n=0}^\infty \frac{H_{n}\binom{2n}{n}}{(n+1)^2\, 2^{2n}}
+2\sum_{n=0}^\infty \frac{\binom{2n}{n}}{(n+1)^2\, 2^{2n}} -\sum_{n=0}^\infty \frac{\binom{2n}{n}}{(n+1)^3\, 2^{2n}}.
\]

Dividing both sides of (\ref{eq:1b}) by $2x,$ and then integrating it yields:
\begin{align*}
&\sum_{n=1}^\infty \frac{\binom{2n}{n}x^n}{2n^{2} \, 2^{2n}}=\int \frac{\log ((1 -\sqrt{1-4x})/2x)}{x} {dx}\\
&=-Li_2\left(\frac{1+\sqrt{1-4x}}{2}\right)-\frac{1}{2} \log^{2}(1+\sqrt{1-4x}) \\
&-\log \left(\frac{1-\sqrt{1-4x}}{2}\right) \log (1+\sqrt{1-4x}) \\
&-\log(-4x) \log \left(\frac{1+\sqrt{1-4x}}{2}\right) +\log(-4x)\log (1+\sqrt{1-4x}) +C.
\end{align*}
This in turn gives us
\begin{equation}
\sum_{n=1}^\infty \frac{\binom{2n} {n}} {2^{2n}}\frac{1}{n^2}=\zeta(2) -2 (\log 2)^2. \label{eq:12}
\end{equation}

Prof Paul Levrie drew my attention to this formula in \cite[p.354, (22)]{srivastava}
\begin{equation}
\sum_{n=1}^\infty \frac{\left(\frac{1}{2}\right)_n}{n \, n!} H_{n-1}=\zeta(2) +2(\log 2)^2  \label{eq:13}
\end{equation}
that involves the Pochhammer symbol $\displaystyle (a)_n =\prod_{k=0}^{n-1}(a+k)$ so that $\left(\frac{1}{2}\right)_n =\displaystyle \frac{\binom{2n}{n}}{n! \, 2^{2n}}.$ The formula can be written as
$\displaystyle \sum_{n=1}^\infty \frac{\binom{2n}{n} \, H_n}{n \, 4^n} -\sum_{n=1}^\infty \frac{\binom{2n}{n}}{n^2 \, 4^n}=\zeta(2) +2(\log 2)^2$
which follows immediately from (\ref{eq:10}) and (\ref{eq:12}).

Now
\[
\sum_{n=1}^\infty \frac{\binom{2n} {n}} {2^{2n}}\frac{1}{n^2}=\sum_{n=0}^\infty \frac{\binom{2n+2}{n+1}}{(n+1)^2 \, 2^{2n+2}}
=\sum_{n=0}^\infty \frac{\binom{2n}{n}}{(n+1)^2 \, 2^{2n}} +\frac{1}{2}\sum_{n=0}^\infty \frac{\binom{2n}{n}}{(n+1)^3 \, 2^{2n+2}}
\]
so using (\ref{eq:6}) and (\ref{eq:12}) we deduce:
\begin{equation}
\sum_{n=0}^\infty \frac{\binom{2n}{n}}{(n+1)^3\, 2^{2n}}=8 -8 \log 2 -\frac{\pi^2}{3} +4 (\log 2)^2. \label{eq:14}
\end{equation}

And by using the sums (\ref{eq:5}), (\ref{eq:6}), (\ref{eq:10}) and (\ref{eq:14}), we obtain:
\begin{equation}
\sum_{n=1}^\infty \frac{H_n \, \binom{2n}{n}}{(n+1)^2\, 2^{2n}} =-\frac{\pi^2}{3} -4 (\log 2)^2 +8 \log 2. \label{eq:15}
\end{equation}

Further, combining (\ref{eq:14}) and (\ref{eq:15}), we get:
\begin{equation}
\sum_{n=1}^\infty \frac{H_n \, \binom{2n}{n}}{n (n+1)^2\, 2^{2n}} =\frac{2\pi^2}{3} +4 (\log 2)^2 -12 \log 2. \label{eq:16}
\end{equation}

Let us replace $x$ by $x^2$ in (\ref{eq:3}) and set $\sqrt{1-4x^2}=y.$ We then take the integral from $y=1$ to $y=0,$ which means taking $x$ from 0 to $1/2.$ Assuming the value of the integral to be \cite[p.173, Table 120, (1)]{haan}
\[
2\int_{0}^1 \log \left(\frac{1+y}{2y}\right) \frac{dy}{\sqrt{1-y^2}}=2\int_{0}^1 \frac{\log (1+u)}{\sqrt{1-u^2}} \, {du}=4G -\pi \log 2,
\]
we obtain:
\begin{equation}
\sum_{n=1}^\infty \frac{H_n \, \binom{2n}{n}}{(2n+1)\, 2^{2n}} =4G -\pi \log 2, \label{eq:17}
\end{equation}
where $G$ is Catalan's constant.

Using a shift of the index, we have:
\begin{align*}
&\sum_{n=1}^\infty \frac{H_n \, \binom{2n}{n}}{(2n-1)^2 \, 2^{2n}} =\sum_{n=0}^\infty \frac{H_{n+1} \, \binom{2n+2}{n+1}}{(2n+1) \, 2^{2n+2}}\\
&=\sum_{n=0}^\infty \frac{H_{n} \binom{2n}{n}}{(2n+1) \, 2^{2n}} -\frac{1}{2}\sum_{n=0}^\infty \frac{H_{n} \binom{2n}{n}}{(n+1) \, 2^{2n}}
+\frac{1}{2}\sum_{n=0}^\infty \frac{\binom{2n}{n}}{(2n+1)(n+1)^2 \, 2^{2n}}.
\end{align*}

By putting $x=1$ in the expansion:
\[
\arcsin x=x+\frac{1}{2}\frac{x^3}{3} +\frac{1\cdot3}{2\cdot4}\frac{x^5}{5} +\frac{1\cdot3\cdot5}{2\cdot4\cdot6}\frac{x^7}{7} +\cdots
\]
we obtain $\displaystyle \sum_{n=0}^\infty\frac{\binom{2n}{n}}{(2n+1) 2^{2n}} =\frac{\pi}{2}.$ We already have two sums $\displaystyle \sum_{n=0}^\infty\frac{\binom{2n}{n}}{(n+1)2^{2n}}$ and $\displaystyle \sum_{n=0}^\infty\frac{\binom{2n}{n}}{(n+1)^2 2^{2n}}$ so that by combining the three sums we get:
\begin{equation}
\sum_{n=0}^\infty \frac{\binom{2n}{n}}{(2n+1)(n+1)^2\, 2^{2n}} =2\pi +4 \log 2 -8. \label{eq:18}
\end{equation}

And using (\ref{eq:17}) and (\ref{eq:18}) we derive this result:
\begin{equation}
\sum_{n=1}^\infty \frac{H_n \, \binom{2n}{n}}{(2n-1)^2\, 2^{2n}} =\pi(1-\log 2) -4(1-G). \label{eq:19}
\end{equation}
where $G$ is Catalan's constant.

Our integral also gives:
\begin{align*}
4\int_{1/3}^{1} \log \left(\frac{1+y}{2y}\right) \, \frac{dy}{1-y^2} &=2\,{\rm Li}_2\left(-\frac{1}{3}\right) +4\,{\rm Li}_2\left(\frac{2}{3}\right)
-\frac{\pi^2}{6}\\ & -(\log 2)^2 +3(\log 3)^2 -4 \, \log 2 \log 3,
\end{align*}
which yields another beautiful formula:
\begin{equation}
\sum_{n=1}^\infty \frac{2^n \, H_n \, \binom{2n}{n}}{n \, 3^{2n}}=\frac{\pi^2}{6} -(\log 2)^2. \label{eq:20}
\end{equation}

Morris\cite[p.781]{morris} notes that $6 {\rm Li}_2(3) -3{\rm Li}_2(-3) = 2\pi^2$ using which and various relations from \cite[p.283]{lewin}, we found: \begin{equation}
2 \, {\rm Li}_2\left(\frac{1}{3}\right) -{\rm Li}_2\left(-\frac{1}{3}\right)=\frac{\pi^2}{6} -\frac{1}{2}\left(\log 3\right)^2 \label{eq:21}
\end{equation}
which we used in the derivation of (\ref{eq:20}). Also refer to \cite[p.155, (2.3)]{kirilov1} \cite[p.89, 6(i)]{kirilov2}.

Further, our integral, in conjunction with the relations \cite[p.283, (7)\&(13)]{lewin}: $\displaystyle {\rm Li}_2(x) +{\rm Li}_2(1-x) =\frac{\pi^2}{6} -\log (x) \log (1-x)$ and $\displaystyle {\rm Li}_2(x) +{\rm Li}_2(-x) =\frac{1}{2}{\rm Li}_2(x),$ yields:
\[
4\int_{1/2}^{1} \log \left(\frac{1+y}{2y}\right) \, \frac{dy}{1-y^2}=\sum_{n=1}^\infty \frac{3^n \, H_n \, \binom{2n}{n}}{n \, 2^{4n}}
={\rm Li}_2\left(\frac{3}{4}\right) +2 (\log 2)^2 -(\log 3)^2
\]
while the relation \cite[p.283, (12)]{lewin} $\displaystyle {\rm Li}_2\left(\frac{1}{1+x}\right) -{\rm Li}_2(-x) =\frac{\pi^2}{6} -\frac{1}{2}\log (1+x) \log \left (\frac{1+x}{x^2}\right) \; x>0$ transforms it into:
\[
\sum_{n=1}^\infty \frac{3^n \, H_n \, \binom{2n}{n}}{n \, 2^{4n}}
=\frac{\pi^2}{6} +{\rm Li}_2\left(-\frac{1}{3}\right) -\frac{1}{2}\left(\log 3\right)^2
\]
and using (\ref{eq:21}) we obtained this lovely result:
\begin{equation}
\sum_{n=1}^\infty \frac{3^n \, H_n \, \binom{2n}{n}}{n \, 2^{4n}}=2 \, {\rm Li}_2\left(\frac{1}{3}\right). \label{eq:22}
\end{equation}
The R.H.S. can also be written as: $\displaystyle \frac{1}{3}{\rm Li}_2\left(\frac{1}{9}\right) +\frac{\pi^2}{9} -\frac{1}{3}\left(\log 3\right)^2.$

\noindent{\textit{Concluding remarks}}: We have given a host of interesting sums here. The reader may try to compute the sums: $\displaystyle \sum_{n=1}^\infty \frac{H_n \, \binom{2n}{n}}{n^2 \, 2^{2n}}$ and $\displaystyle \sum_{n=1}^\infty \frac{H_n \, \binom{2n}{n}}{(2n+1)^2 \, 2^{2n}}$ which we couldn't.

\noindent {\bf Acknowledgement:} The author is thankful to Prof K. N. Boyadziev for suggesting explicit use of generating functions, and to Prof P. Levrie for referring to a formula in \cite {srivastava}. Ming Yean detected few errors in the version \url{https://arxiv.org/abs/1806.03998v1}. The author is grateful to the learned referee who too pointed to the same and also made some useful suggestions regarding language.

\end{document}